\documentclass[12pt]{amsart}
\usepackage{graphicx}

\usepackage{amssymb, microtype}

\newcommand{\color}[1]{{}}

\title[The discovery of hyperbolic structures]
{A personal account of the discovery of hyperbolic structures on some knot complements}

\author{Robert Riley}

\begin{document}

\maketitle

\vspace{-25pt}

\noindent  Mathematical Sciences, Binghamton University, Binghamton,
NY, USA 

\vspace{15pt}

{MSC classification (2010): 
\newline Primary: 57M50, 57M25, 
\newline Secondary 30F40,
01A60} 

\vspace{10pt}

Keywords:
Hyperbolic structures,
Knot complements

Abstract: I give my view of the early history of the discovery of
hyperbolic structures on knot complements from my early work on
representations of knot groups into matrix groups to my meeting with
William Thurston in 1976\footnote{This article was written by Robert
Riley about ten years before his death in 2000 and never submitted
for publication.  An explantion of why it is being published now and
some information about Riley and this article is given in \cite{BJS}
which accompanies this article in this issue of the journal.}

\section{Introduction}

I discovered, quite unexpectedly, the phenomenon of hyperbolic structure on
three knot complements early in 1974, and managed to get two papers on the
topic published in 1975.  At some moment between the dates of publication of
these papers, William Thurston independently discovered the phenomenon and
ran away with the idea.  In late June or early July 1976 he learned of my
work, and so when we met later in July he immediately told me that he had
been trying for about a year to prove the hyperbolization conjecture for
Haken 3--manifolds.

Colin Adams published a semipopular account of knot theory in ``The
Knot Book'' ~\cite{Ad}, and a copy of this came into my hands
recently.  On page 119 he gives an account of the hyperbolic
structure discovery which is just plain wrong.\footnote{Riley refers
here to the first edition of {\itshape The Knot Book}, published in
1994.  In \S5.3 of the 2004 edition, published by the AMS, there is
a concise, corrected account of this discovery, together with an
excellent elementary introduction to hyperbolic knots.} He does get
the names of the two people concerned and the priority right, but
nothing else.  The present paper is an attempt to set the record
straight.  I shall relate what I did, why, and when.  There will be
too much detail about small matters, but this will convey the spirit
of my projects.  Indeed, I think my old papers were very open about
my project, and a close look at them and their dates of submission
should have made the present history unnecessary.  Furthermore, Bill
Thurston's account of my work in~\cite{Th1} is entirely fair, except
for being too generous about my influence on his thinking.

So below I give the history of my project from its beginning to the moment I
met Professor Thurston.  The story is told as I saw it, and the emphasis is
on motivation and dates.  Many dates are only approximate because most
entries in my notebooks are undated, but the uncertainties are never more
than about a month.  I include an intermediate example, worked out between
the discoveries of the hyperbolic structures for the figure--eight knot
$(4_1)$ and for $5_2$.  This example ought not be on the main line of
development, but in fact it was, and it served to undermine my initial
expectation that the figure--eight is the only knot which could possibly be
hyperbolic.  I close with some comments on the early work of H.~Gieseking
{\color{red}and} Max Dehn, and on the article~\cite{Th3} of W.~Thurston.  


\section{The early years}

On settling in Amsterdam in October 1966 I wrote off to virtually everyone
publishing in knot theory for their reprints and preprints.  I recall with
gratitude that R.H.~Fox and H.~Seifert were especially generous.  An
unassuming little paper by Fox~\cite{Fo} written in Utrecht some 20 miles away,
took my fancy.  Here Fox advertised the notion of longitude in a knot group
by using it, together with representations on the alternating group $A_5$,
to distinguish the square and granny knots.  I was intrigued by the success
of $A_5$, and took the first steps toward writing out explicit procedures to
find all $A_5$--representations of a knot group in 1967--68.  When I got my
first temporary appointment at Southampton (England) in 1968 this became my
main project, with results summarized in~\cite{Ri1, Ri2}.  So by 1970 I was after the parabolic representations ($p$--reps) of a knot group, initially because they
were easier to manage than the general non--abelian representations
($nab$--reps).  The 2--bridge case is especially tractable, because the
representations are governed by a simple polynomial whose rule of formation
is easily programmed in Fortran.  This tractability extends to all
$r$--bridge knots which are symmetric about an $r$--fold axis of rotation that
cyclically permutes the bridges, but most knots of bridge number $> 2$ are
not so symmetric.  The explicit algebraic description of the equivalence
classes of $p$--reps of an unsymmetric knot is so difficult that only a few
examples have been worked explicitly, and I have found the full curve of all
$nab$--representations of only one unsymmetric 3--bridge knot, $8_{20}$.  Around 1971 I wrote some primitive Fortran programs to find the $p$--reps of a few 3--bridge knots and used the output to discover the commuting trick of~\cite{Ri2},~\cite[II]{Ri6}, but at the time this topic was mainly pure frustration. 

In 1971 a plea for help from me was passed on to Professor
G.E.~Collins, the instigator of the SAC--1 file of Fortran routines
for doing the kind of algebraic calculations I needed.  He sent me a
pile of very poorly printed manuals containing the program listings,
lots of errata slips, and the advice that the 24 bit word size of
the Southampton {\color{red}U}niversity computer would require some
doubly recursive programming in assembly language.  (The reference
count field in a SAC atom would be too small without this recursion,
and hence impose a strict limit on the allowed complexity of
calculations).  He also mentioned that I would need to get someone
to punch up the 6000 or so cards of the 1971 SAC.  Well, that
someone had to be me, but fortunately only some 4000 cards, plus the
assembly language parts, were needed for my application.  It took
about eight months to do all this, and I never did get the double
recursion for the reference count right.  So my more ambitious
calculations were killed as soon as the reference count tried to
reach 128, but I still managed to do most of what I wanted.  By 1
October 1972, the day my fourth temporary appointment at Southampton
ceased, I had done the elimination--of--variables part of the
solving for an algebraic description of the set of $p$--reps for
several 3--bridge knots, including $9_{35}$.  Each SAC run required
several hours of CPU time, and could not have been attempted during
term time.  Perhaps some distorted memory of this story is the
source of the ``immense computer program that was designed to
attempt to show that some knots are hyperbolic" bit in Adams'
account.  In fact, the PNCRE package which does just this was
developed from 1976, and it was always fast enough for term time,
even during the day on a grossly overloaded 1960's computer.


\section{The preparation}

In October 1972 I had a large pile of SAC output which needed more computer
analysis to become meaningful, and no prospect of further employment.  So I
spent the next three months walking the Pennine Way and walking in Wales
until the prospect of a six month appointment in Strasbourg opened up.
While I was walking in the Vosges this materialized, and I was able to
complete the algebraic description of the equivalence classes of $p$--reps
for several knots, including $9_{35}$, cf.~[11].  (I recall a puzzling
difficulty with $9_{32}$ that was explained a decade later as the
consequence of dropping the deck of data cards, perhaps in 1971, and
reassembling it almost exactly right).

The knot $9_{35}$ has a large symmetry group (dihedral of order 12, [11]),
and also an unusually large number of algebraic equivalence classes of
$p$--reps, facts which I believe are related.  The SAC calculations had given
me a polynomial $p(x) \in \Bbb Z[x]$ of degree 25 which I had to factor as
the first step.  When one has no symbolic manipulation package available
this is done by finding the roots of $p(x) = 0$ and examining them for
clues.  The polynomial $p(x)$ (and its relative for $9_{48}$ which was even
worse) defeated several commercially produced root--finding routines, but a
final resort routine succeeded, sort of, and I was able to infer factors 
$$p_1 = 1 + x,\ \ p_2 = 1 + 2x + 7x^2 + 5x^3 + x^4,\ \ p_3 = \cdots,$$
and soon
$$p(x) = (1+x)^{10}p_2(x)^2\cdots.$$
Only the cubic factor remained unguessed, and of course it turned out to be
the one giving the hyperbolic structure four years later.  Each factor
$p_k(x)$ of $p(x)$ had to be tested to see if it gave an equivalence class
of $p$--reps or was spurious, and I expected $1 + x$ to be spurious.  To my
surprise it gave $p$--reps on 
$$G_i = \left\langle \bmatrix 1 &1\\ 0 &1\endbmatrix,\  \bmatrix 1 &0\\ -1
&1\endbmatrix,\  A_i\bmatrix 1 &0\\ -1 &1\endbmatrix A^{-1}_i\right\rangle
\subset SL_2(\Bbb Z[i]),\  A_i = \bmatrix 1 &i\\ 0 &1\endbmatrix,$$
where $i = \sqrt{-1}$.  This was in June 1973, and I probably did not
understand what a Kleinian group is at the time, but I could see $G_i$ is
discrete and wondered what its presentation was.  Also, as I watched the
printout emerge from the line printer I guessed that these $p$--reps must be
an instance of an undiscovered theorem, and the same evening stated and
proved the theorem.  (Writing it up for publication is taking longer.  In
December 1991 I used Maple to extend the theorem to algebraic varieties of
$nab$--reps and add some new material.  In 1993 I told Tomotada Ohtsuki
about this, giving no detail, and he promptly found a better proof and more
new material.  I hope to proceed to a joint paper soon.)

After the summer vacation of 1973 when I returned to Southampton,
the professors of the mathematics department granted me the use of
an office and all university facilities, except the computer which
was {\it heavily} overloaded.  By then I had learned by some osmosis
what a Kleinian group is and read Maskit's paper~\cite{Ma} on
Poincar{\'e}'s Theorem on Fundamental Polyhedra.  This made progress
on $G_i$ above possible, and I soon had its presentation.  (I also
found that Fricke and Klein had considered $G_i$, or something very
like it, cf.~Fig.~151 on page 452 of~\cite{FK}{\color{red}.)}
Success with $G_i$ led to success with the image $\pi K\theta$ of a
$p$--rep of the figure--eight knot group in November 1973.  Recall
that $$\pi K\theta = \left\langle \bmatrix 1 &1\\ 0 &1\endbmatrix,
\bmatrix 1 &0\\ -\omega &1\endbmatrix\right\rangle,\ \ \omega =
\frac{-1+\sqrt{-3}}{2},$$ so the group is obviously discrete and
only its presentation was in doubt.  I remember my surprise at
finding this $p$--rep is faithful.  The first version of my
account~\cite{Ri3} of this was received by the Editors on 30
November 1973, and it didn't mention the orbit space $\Bbb H^3/\pi
K\theta$ because I had not even thought of it.

Why not?!  Well, the result was perhaps a fortnight old, and I didn't have a
premonition of hyperbolic structure on knot complements.  Years later I
learned that it had not only been thought of, but attempted and discussed
privately by the Kleinian groupies since 1968.  Nothing had been written and
none of this had reached me.  The key to seeing that the orbit space of $\pi
K\theta$ had to be the figure--eight complement was seeing the peripheral
torus in the orbit space.  This torus occurs as the image of Euclidian plane
$\Pi(t) = \{(z,t) : z\in \Bbb C\}\subset \Bbb H^3$ for any $t > 1$.  In my
diagram $\Pi(t)$ meets the fundamental domain not in a parallelogram but in
a zigzag shape (four hexagonal discs), and perhaps the zigzag temporarily
prevented me from seeing the torus.  This is silly, because the stabilizer
of the torus is the free abelian group $(\pi K\theta)_\infty$ generated by
$z\mapsto z + 1$, $z\mapsto z + 2\sqrt{-3}$, and $(\pi K\theta)_\infty$ has
to be considered explicitly during the verification that Poincar{\'e}'s
theorem applies to my supposed Ford fundamental domain.  But silly or not,
it took perhaps seven weeks, till January 1974, for me to see the torus.
Verification that $\Bbb H^3/\pi K\theta = S^3 - \text{fig--eight}$ took
perhaps a day, and consisted of looking at my reprint of Waldhausen's paper~\cite{Wa}.  It seems unfortunate that this was too easy, and that I should have
been forced to develop a direct geometrical argument, but once the pressure
was off I didn't want to do it.  I expect that a direct geometrical
construction works for all non--torus two bridge knots, and that it would
prove the conjectures of~\cite[\S4]{Ri7}, so the matter will not be a waste of
effort.

The figure--eight discovery was not decisive for me as it was for Thurston.
I expected that Shimizu's lemma, viz.
$$\left\langle \bmatrix 1 &1\\ 0 &1\endbmatrix, \bmatrix a &b\\ c
&d\endbmatrix\right\rangle\ \text{ is not discrete when }\ ad -bc = 1,\ \ 0
< |c| < 1,$$
would preclude the discreteness of the images $\pi K\theta$ of the
potentially faithful $p$--reps $\theta$ for all other knots.  (In particular,
I predicted Alan Reid's theorem~\cite{Re} that the figure--eight is the only
arithmetic hyperbolic knot).  However, by the time I mailed off the revised
version of~\cite{Ri3} that was actually accepted I had recognized the true
situation, but, I suppose out of laziness, I didn't revise~\cite{Ri3} again to make
an announcement.

R.H.~Fox died within a few days of the figure--eight discovery.


\section{The intermediate example}

I now had a beautiful discovery, and a certain fear of testing
whether something similar was true for the obvious next case, the
knot $5_2$ of two--bridge types $(7,3)$, $(7,5)$.  Instead of going
for $5_2$ directly I temporized by taking up a different kind of
example, the groups $\pi K\theta$ associated to a cubic factor
$f(u)$ of the $p$--rep polynomial for the knot $8_{11}$ of
two{\color{red}--}bridge types $(27,17)$, $(27,19)$.  To give an
account of this we need to recall the basics of
two{\color{red}--}bridge knot groups and their $p$--reps.

A two{\color{red}--}bridge knot normal form corresponds to a pair
$(\alpha,\beta)$ of integers, where $\alpha > 1$ is odd, $\beta$ is
odd, $gcd(\alpha,\beta) = 1$, and $-\alpha < \beta < \alpha$.  The
knot group $\pi K$ for $(\alpha,\beta)$ depends not on $\beta$
itself but on $|\beta|$, so we may as well assume $0 < \beta <
\alpha$.  Then $$\pi K = |x_1,x_2 : wx_1 = x_2w|, \ \ w =
x_1^{\epsilon_1}x^{\epsilon_2}_2\cdots
x^{\epsilon_{\alpha-1}}_2,\eqno(3.1)$$ where $\epsilon_j =
\epsilon_{\alpha-j} = \pm 1$, and the exponent sequence $\overset
\rightharpoonup \epsilon = (\epsilon_1,\cdots,\epsilon_{\alpha-1})$
is determined by a simple rule, cf{\color{red}.}~[7, 12].  A
longitude $\gamma_1$ in the peripheral subgroup $\langle
x_1,\gamma_1\rangle$ of $x_1$ is a certain word $\tilde w^{-1}
wx_1^{2\sigma}$ on $x_1,x_2$.  A normalized $p$--rep $\theta =
\theta(\omega)$ of $\pi K$ is a homomorphism such that $$x_1\theta =
A = \bmatrix 1 &1\\ 0 &1\endbmatrix,\ \ x_2\theta = B = B_\omega =
\bmatrix 1 &0\\ -\omega &1\endbmatrix,\eqno(3.2)$$ where $\omega \in
\Bbb C$.  Indeed, $\omega$ is a root of the $p$--rep polynomial
$\Lambda[u] \in \Bbb Z[u]$ which may be reducible but which has no
repeated roots.  Then the longitude entry $g(\theta)$ or $g(\omega)$
for $\theta(\omega)$ is found by $$\gamma_1\theta = \bmatrix -1
&g(\theta)\\ 0 &-1\endbmatrix,$$ and is readily computable once
$\omega$ is known.  To factor $\Lambda(u)$ without a system like
SAC, Macsyma, or Maple but when a polynomial root finding package is
available{\color{red},} find the roots and list the pairs
$(\omega,g(\omega))$.  Factors stand out as having pairs where
$g(\omega)$ evidently belongs to a proper subfield of $\Bbb
Q(\omega)$.  In the case of $8_{11}$ we found the factor $$f(u) = -1
+ u(1+u)^2$$ by $g(\theta) = -6$ for its roots.  The roots of $f(u)$
are $$\omega_1 = -1.23278 + 0.79255i, \ \ \omega_2 = \bar\omega_1,\
\ \omega_3 = 0.46557,\eqno(3.3)$$ (rounded to 5 decimal accuracy).
Today this factor is explained as an instance of Theorem B
of~\cite{Ri7} and it clearly had something to do with the discovery
of the theorem.  I had $f(u)$ by 1971.

By February 1974 my worries about the figure--eight knot brought me to
consider the group $\Gamma = \langle A,B\rangle$, $B= B_\omega$, where
$\omega$ is the $\omega_1$ of (3.3).
I simply went for a Ford domain $\mathcal D$ of $\Gamma$ using graph paper, compass and ruler, and the first programmable calculator available at Southampton.
(That would have cost about
two months gross salary if I had still been employed).  It didn't take long
to get the diagram of Fig.~1, and when the time came to think about proof
the closing trick and angle sum trick of~\cite{Ri5} came to mind automatically.
As far as I know this group $\Gamma$ is the first group proved discrete by
Poincar{\'e}'s theorem where these tricks are necessary.  Perhaps the first
people to wonder about using Poincar{\'e}'s theorem for computation with
potentially discrete groups didn't see these simple tricks in advance,
didn't have a specific example they really needed, and shied away from
getting too involved.

\begin{figure}
\begin{center}
\includegraphics[scale=0.7]{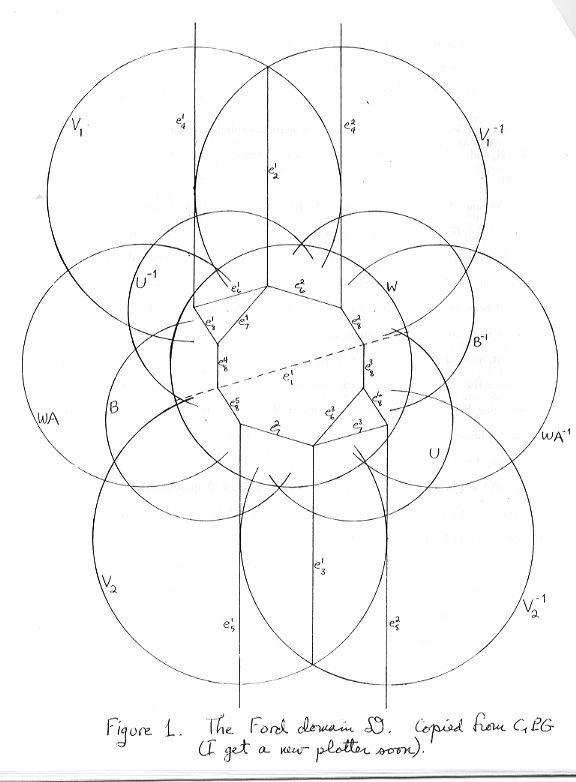}
\end{center}
\end{figure}

We give a little more detail on $\Gamma$ and its Ford domain
$\mathcal D$ illustrated in Fig.~1.  This is taken from an
unpublished paper CPG, written in late 1974 and early 1975, doing
all the discrete non{\color{red}--}Fuchsian cases where the group
$\pi K\theta$ corresponds to a root of a cubic polynomial,
viz.~$5_2$, $7_4$, and $8_{11}$.  The case $5_2$ is worked in
\cite{Ri6}, and $7_4$ is similar to but easier than
$7_7${\color{red},} also worked in~\cite{Ri6} but much easier.

Let $\pi K$ be the group of $(27,17)$ presented as in (3.1), so $\Gamma =
\pi K\theta$ as in (3.2).  We have words
$$u := x^{-1}_1x_2x_1,\ \ v_1 := u x^{-1}_2x_1x^{-1}_2,\ \ w_1 :=
v_1x^{-1}_1v^{-2}_1x^{-1}_2v_1x_1^{-1}.$$
The word $w$ of (3.1) is $w_1x_1$, so $w_1x_1 = x_2w_1$ holds in $\pi K$.
These words $u$, $v_1$ were found by straightforward search of subsegments
of $w$ to correspond to spheres carrying sides of tentative Ford domains.
The search for the sides of a fundamental domain has to be guided by some
principle, since a Cantorian exhaustion is too slow, and segments of $w$
worked well, both here and later for all two bridge knots.

We found easily that the elements 
$$A = x_1\theta,\ \ U = u\theta,\ \ V_1 = v_1\theta,\ \ W_1 = w_1\theta,\ \
V_2 = U^{-1}W_1$$
seem to be the side pairing transformations of the tentative Ford domain
$\mathcal D$ of Fig.~1.
Thus we read off from Fig.~1 a proposed presentation
for $\Gamma$:  generators $A$, $U$, $V_1$, $V_2$, $W_1$.  relations
$$\begin{array}{ll}
&\quad W^2_1 = V^3_1 = V^3_2 = (A^{-1}V_1)^2 = (A^{-1}V_2)^2 = E,\\
&V_1=  W_1U^{-1},\ \ V_2 = U^{-1}W_1,\ \ U = A^{-1}W_1AW_1A^{-1}.\\
\end{array}
$$
To use the closing trick and angle sum tricks of [10] it is necessary to
verify directly that these relations hold in $\Gamma$.  For this it helps to
see copies of the modular group $SL_2(\Bbb Z)$ in $\Gamma$.  Let
$$A_* := \bmatrix 1 &u + u^2\\ 0 &1\endbmatrix,$$
then
$$V_1 \equiv A^{-1}_*\bmatrix 0 &-1\\ 1 &1\endbmatrix A_*,\ \ V_2 \equiv
A_*\bmatrix 1 &-1\\ 1 &0\endbmatrix A^{-1}_*\ \text{(mod }f(u)).$$
So $\langle A,V_1\rangle$ and $\langle A,V_2\rangle$ are conjugate to
$SL_2(\Bbb Z)$ in $SL_2(\Bbb Z[\omega])$.  All the proposed relations now
can be verified by straightforward computation in $SL_2(\Bbb Z[u])$ modulo
$f(u)$.  Then the arguments of [10] show that $\Gamma$ is discrete, $\mathcal D$
is a fundamental domain for it, and that these relations present the group.
This made a good confidence--building exercise for me, and might do the
same for other people.  Note that this $\mathcal D$ is simpler than the Ford
domain for $5_2$ discussed in~\cite{Ri6}, so $\Gamma$ really is an intermediate
example.  


\section{Completion of the discovery}

This procrastination had now given me a bigger worry which can be put thus:
Why should the Great Lord have performed a unique miracle to make $\Gamma$
discrete, for no visible reason at all?!  The answer is compelling:  He
didn't!  If $\Gamma$ is discrete then many other groups have to be discrete,
in direct defiance of Shimizu's lemma, and, since each case of discreteness
requires a good reason, there must be general theorems explaining this
discreteness.  It is a little ironic that this prediction was amply
vindicated for (suitable) 3--manifold groups, but, at this writing, the
general theorem explaining the discreteness of $\Gamma$ has not been stated,
let alone proved.  

During a few weeks further procrastination the above considerations
compelled me to predict that a knot in $S^3$ is hyperbolic unless it clearly
was not.  Early in March 1974, I think, I finally went to work on $5_2$ and
in a few hours had confirmed my prediction.  This completed the essential
part of my discovery, and all later cases, such as $7_4$ and several links,
were just routine examples at most illustrating matters of secondary
importance, such as the symmetries of a knot.  In fact, for a while I was
confused by the symmetries and thought that a too--rich symmetry group would
preclude the hyperbolic structure, but I eventually found my mistake.  So by
late 1974 I had gotten it right:  a knot is hyperbolic unless its group
contains a noncyclic abelian subgroup which is not peripheral.  Making bold
sweeping conjectures is unnatural for me, and I didn't venture to predict
anything about arbitrary 3--manifolds.  I suppose that I might have
predicted which 3--manifolds were hyperbolic had someone pressed me on the
issue in conversation, but I was too isolated and unknown for that to
happen.  The locals at Southampton were rather cool about the whole project,
except for David Singerman.  He liked it enough to propose that we try to
get the Science Research Council (of Great Britain) to support me on a
hyperbolic project at Southampton University while I got my Ph.D. and looked
for a permanent job.  His plan was to time the submission of the proposal so
that the referee would be at the summer 1975 conference on Kleinian groups
at Cambridge where I would publicize hyperbolic structure.  Whether or not
the plan worked, the Kleinian groupies liked my examples, especially because
these examples pointed up the importance of their own work.  The SRC did fund
the project generously, ultimately for four years 1976--1979.

The first two years of the project were devoted to the development
of the system PNCRE~\cite{Ri5}, a file of Fortran subroutines to
compute with explicit subgroups of $SL_2(\Bbb C)$.  PNCRE was not
easy to develop and its first output came early in 1977.  Meanwhile,
about March 1976, a colleague gave me a preprint of Thurston's
lecture~\cite{Th1} on fo{\color{red}l}iations of surfaces.  This was
the first I heard of him, and I recall that on reading it I became
certain that he and I would never share any common mathematical
interest.  In late June 1976 a friend drove me up to the University
of Warwick to hear a lecture by J.~Milnor on topics like
Sarkovskii's theorem.  Directly he was finished I very nervously
(read: scared stiff) introduced myself to him and told him about
examples of hyperbolic knots/links.  He was interested, and asked a
number of direct questions, so that in a minute he understood the
status of my project (examples only).  I did not guess that he
already knew something about the matter.  I was so scared that when
he asked me to repeat my name I simply ran away.  But perhaps even
before we got back to Southampton that evening, Milnor had asked the
locals who in Britain was interested in hyperbolic structure on knot
complements, and directly afterwards Thurston had his hands on my
two papers.  If not, he did when I sent my papers to Milnor the next
week (early July 1976).

Later that month I was invited, to put it mildly, to spend a week in David
Fowler's home in Warwick.  His wife is French, and she felt that that year
she simply had to bring the children to France to meet their relatives.  She
naturally had the house and garden filled with beautiful plants which need
constant watering.  The summer of '76 was a famous drought in which the
water shortage was so severe that the only legal water for plants was used
bath water.  Hence the urgent need to have the Fowler's home occupied every
night, and David Rand, who had been a student at Southampton and was taking
up a lectureship at Warwick, put me down for one week.

On my arrival in the common room of the Warwick Mathematics Department,
David Epstein sprang up and asked me who I was.  He had seen my face on
numerous occasions over the years, most recently when I sat directly behind
him in Milnor's lecture, and wanted to know.  On hearing my name, a tall man
sprawled over three chairs sprang up.  He said he was Bill Thurston, that he
wanted to meet me, and that for about a year he had been working on a
general conjecture which included everything I was doing.  The shock was
immense.  I am afraid that I react badly to surprises, and I became quite
unpleasant for the rest of the week.  Fortunately Bill didn't hold it
against me later.  His later statement (page 177 of~\cite{Th3}),

``{\color{red}$\ldots$}; and I have not actively or effectively
promoted the field or the careers of the excellent people in it.''

was either not written with me in mind or he judges me not to
satisfy the qualification. He certainly did advance my career actively:
strong letters of recommendation, several thousands of dollars from his
Waterman Fellowship, inclusion in the 1980--81 Thurston--Sullivan NSF
project at Boulder, and a trip to Binghamton at my request.  I owe
everything to the people who have so generously supported me over the years
when I needed help most:  H.B.~Griffiths, David Singerman, and Bill
Thurston, and I am deeply grateful to them all.


\section{Hindsight}

The question is:  Why did the explicit discovery of hyperbolic structure on
at least some knot complements wait until 1974?  Wilhelm Magnus told us that
H.~Gieseking, in a thesis written in 1912 under the direction of Max Dehn,
considered a group $G_1$ of hyperbolic isometries of a ball model $\mathcal B^3$
of $\Bbb H^3$ and certain of its subgroups.  The fundamental domain for
$G_1$ is a regular ideal tetrahedron $T_1$, and $G_1$
contains orientation reversing elements.  Gieseking considered the
orientation preserving subgroup $G_2$ of index two whose fundamental domain
is two tetrahedra glued together along a face, without recognizing that
$G_2$ is isomorphic to the figure--eight knot group and its orbit space is
the figure--eight complement.  Magnus told me that Dehn considered these
groups only as exercises in geometric symmetry:  the geometric description
of $G_1$, $T_1$ is so simple that Poincar{\'e}'s theorem simply has to apply
directly.  If Dehn had known that the figure--eight knot was involved he
certainly would have had Gieseking publish, he would have given the matter
the greatest publicity, and the development of 3--manifold theory would have
gone very differently.  So why did they not recognize the figure--eight
complement?

I propose an answer to this question analogous to my own experience:  I
didn't see the peripheral torus for several weeks but when I did I knew what
I had to have.  I began with the figure--eight and had it in mind.  They
began with an exercise in symmetry and had nothing further in mind.
Furthermore they would have to ask the question:  for $\epsilon > 0$ let
$S_\epsilon$ be the 2--sphere in $\mathcal B^3$ with centre $\overset
\rightharpoonup  0$ and radius $1-\epsilon$.  Then $G_1$ maps
$S_\epsilon$ to itself, so what is the orbit space $S_\epsilon/G_1$?  With
hindsight the answer is obvious:  a Klein bottle.  Dehn would have answered
this question easily once it had been raised, and I feel certain the Klein
bottle would have disturbed him deeply.  The result would have burned within
him until he was driven to get to the bottom of the matter, and somehow he
would have found the figure--eight.  They had about two years to do this
before the Great War of 1914--18 swept Gieseking to his doom.  Dehn's good
students were probably all destroyed, and most likely Dehn was so distraught
at their loss that he couldn't bear to think about his joint projects with
them any longer.

As far as I know, the next time critical examples of hyperbolic structure
on 3--manifolds should have been found was in the late 1950's, during the
period of euphoria caused by Papakyrikopoulos' breakthrough with his proofs
of Dehn's lemma etc.  The topic was definitely thought of, but nothing
happened, perhaps because the man concerned did not have anything specific
to work on, and he certainly had a lot of other important projects to
pursue.  In 1968 a Kleinian groupie wondered whether a knot complement could
be hyperbolic, and chose as example to test this idea the trefoil knot.  He
soon found it didn't work and was discouraged.  (Actually, the trefoil
complement does carry hyperbolic orbifold structures of infinite volume, but
nobody wanted that).  In the early 1970's he actually visited Southampton
University and met me, but somehow the crucial topic didn't come up in the
discussion.  If it had I would have put him onto the figure--eight and even
given him the exact matrices to use.  I could not have done the calculation
with Poincar{\'e}'s theorem at the time (he could), but I did have
Waldhausen's paper to help with identifying the orbit space.

I would like to close by quoting a paragraph from page 175 of Thurston's
essay~\cite{Th3}.

``Neither the geometrization conjecture nor its proof for Haken
manifolds was in the path of any group of mathematicians at the time --- it
went against the trends in topology for the preceding 30 years, and it took
people by surprise.  To most topologists at the time, hyperbolic geometry
was an arcane side branch of mathematics, although there were other groups
of mathematicians such as differential geometers who did understand it from
certain points of view.  It took topologists a while just to understand what
the geometrization conjecture meant, what it was good for, and why it was
relevant.''

Well, this is not quite right.  For one thing, it is too strongly
put.  When I met them in 1975 the Kleinian groupies had been knowledgeable
about the hyperbolization conjecture for Haken manifolds for at least a
couple of years, but they saw it as too much for themselves.  For another,
what really took people aback was the speed with which the task was completed
(excepting the write--up).  Thurston simply didn't give anyone starting from
my examples the time to get involved.  And few serious mathematicians would
look at one modest example of something pretty and
immediately formulate the most sweeping conjecture for 3--manifolds which
could possibly be true, and then plunge in.  Thurston's success at doing
this is his own personal triumph, and not a closing out of a golden
opportunity that the rest of us were fool enough to lose.

I had thought of saying something about the history of Bill Thurston's
thinking about hyperbolic structure in the two years before we met, but I am
afraid to repeat Colin Adam's mistake.  There are rumours that he initially
thought the hyperbolic structure for the figure--eight was impossible,
because of difficulties with the lift of a Seifert surface to $\Bbb H^3$,
and that he discussed these matters with William Jaco at a conference.  The
story continues that when Bill got back to Princeton he found his supposed
contradiction disappear (the lift of the Seifert surface meets the sphere at
infinity in a Peano curve), that this completely reversed his expectations,
and that he first got the figure--eight out of an example of Troels
J{/\!\!\!o}rgensen.  I cannot vouch for any of this.  






\end{document}